\documentclass{gtart}

\def\ifplaintex{\expandafter\ifx\csname documentclass\endcsname\relax}

\def\gtp{{\mathsurround=0pt\it $\cal G\mskip-2mu$eometry \&\ 
$\cal T\!\!$opology $\cal P\!$ublications}}  

\def\recd{{\small Received:\qua\receiveddate\ifx\reviseddate\relax
\else\qquad Revised:\qua\reviseddate\fi\par}} 

\def\lognumber#1{\def\thelognumber{#1}}
\def\volumenumber#1{\def\thevolumenumber{#1}}
\def\volumeyear#1{\def\thevolumeyear{#1}}
\def\papernumber#1{\def\thepapernumber{#1}}
\def\pagenumbers#1#2{\def\startpage{#1}\def\finishpage{#2}}
\def\published#1{\def\publishdate{#1}}

\def\received#1{\def\receiveddate{#1}}
\def\revised#1{\def\reviseddate{#1}}
\def\accepted#1{\def\accepteddate{#1}}

\def\asciiaddress#1{\def\theasciiaddress{#1}}

\let\\\par\let\thelognumber\relax\let\thevolumenumber\relax
\let\thepapernumber\relax\let\thevolumeyear\relax\let\startpage\relax
\let\finishpage\relax\let\publishdate\relax\let\receiveddate\relax
\let\reviseddate\relax\let\accepteddate\relax\let\theasciititle\relax
\let\theasciiauthors\relax\let\theasciiaddress\relax
\let\theasciiabstract\relax

\let\theasciiemail\relax

\ifplaintex
\font\logobig=cmssbx10 scaled 3836
\font\logomed=cmssbx10 scaled 2557
\else
\font\logobig=cmssbx10 scaled 4200
\font\logomed=cmssbx10 scaled 2800
\fi

\long\def\makeagttitle{   
\count0=\startpage
\agt\hfill      
\hbox to 45truept{\vbox to 0pt{\vglue -13truept{\logomed A\kern -.37em{\logobig 
T}\kern -.38em G}\vss}\hss}
\break
{\small Volume \thevolumenumber\ (\thevolumeyear)
\startpage--\finishpage\nl
Published: \publishdate}

\vglue .25truein

{\parskip=0pt\leftskip 0pt plus
1fil\def\\{\par\smallskip}{\Large\bf\thetitle}\par\medskip} \vglue
0.05truein

{\parskip=0pt\leftskip 0pt plus 1fil\def\\{\par}{\sc\theauthors}
\par\medskip}%
 
\vglue 0.03truein 

{\small\leftskip 25truept\rightskip 25truept{\bf Abstract}\stdspace\theabstract

{\bf AMS Classification}\stdspace\theprimaryclass
\ifx\thesecondaryclass\relax\else; \thesecondaryclass\fi\par
{\bf Keywords}\stdspace \thekeywords\par}\vglue 7truept

}   

\ifplaintex
\font\phead=cmsl9 scaled 950
\font\pnum=cmbx10 scaled 913
\font\pfoot=cmsl9 scaled 950
\headline{\vbox to 0pt{\vskip -4.5mm\line{\small\phead\ifnum
\count0=\startpage ISSN 1472-2739 (on-line) 1472-2747 (printed)
\hfill {\pnum\folio}\else\ifodd\count0\def\\{ }%
\ifx\theshorttitle\relax\thetitle\else\theshorttitle\fi\hfill{\pnum\folio}
\else\def\\{ and }{\pnum\folio}\hfill\ifx\theshortauthors\relax\theauthors
\else\theshortauthors\fi\fi\fi}\vss}}
\footline{\vbox to 0pt{\vglue 0mm\line{\small\pfoot\ifnum\count0=\startpage
\copyright\ \gtp\hfill\else
\agt, Volume \thevolumenumber\ (\thevolumeyear)\hfill\fi}\vss}}
\else
\font=cmsl9 scaled 1050
\font\lnum=cmbx10 
\font=cmsl9 scaled 1050
\makeatletter
\def\@oddhead{{\small\lhead\ifnum\count0=\startpage ISSN 1472-2739 
(on-line) 1472-2747 (printed)\hfill {\lnum\number\count0}\else\ifodd\count0
\def\\{ }\ifx\theshorttitle\relax \thetitle \else\theshorttitle\fi\hfill
{\lnum\number\count0}\else\def\\{ and }{\lnum\number\count0}
\hfill\ifx\theshortauthors\relax 
\theauthors\else\theshortauthors\fi\fi\fi}}\def\@evenhead{\@oddhead}
\def\@oddfoot{\small\lfoot\ifnum\count0=\startpage\copyright\ \gtp\hfill\else
\agt, Volume \thevolumenumber\ (\thevolumeyear)\hfill\fi}
\def\@evenfoot{\@oddfoot}
\makeatother
\fi
\let\maketitlepage\makeagttitle
\let\makeshorttitle\maketitlepage
\let\maketitle\maketitlepage

\newwrite\gtoutfile
\long\gdef\makeheadfile{  
{\def\\{, }\def\s{ }
\immediate\openout\gtoutfile head.xxx
\immediate\write\gtoutfile{To: math@arxiv.org}
\immediate\write\gtoutfile{Subject: put OR rep NNNNN:ppppp}
\immediate\write\gtoutfile{--text follows this line--}
\immediate\write\gtoutfile{Proxy-for: \ifx\theasciiauthors\relax
\theauthors\else\theasciiauthors\fi\s<\ifx\theasciiemail\relax\theemail\else\theasciiemail\fi>}
\immediate\write\gtoutfile{\noexpand\\}
\immediate\write\gtoutfile{Authors: \ifx\theasciiauthors\relax
\theauthors\else\theasciiauthors\fi}
{\def\\{ }\immediate\write\gtoutfile{Title: \ifx\theasciititle\relax
\thetitle\else\theasciititle\fi}}
\immediate\write\gtoutfile{Subj-class: GT or SG, GR etc}
\immediate\write\gtoutfile{MSC-class: \theprimaryclass\ifx\thesecondaryclass\relax\else, \thesecondaryclass\fi}
\immediate\write\gtoutfile{Journal-ref: Algebr. Geom. Topol. \thevolumenumber\s
(\thevolumeyear) \startpage-\finishpage}
\immediate\write\gtoutfile{Comments: Published by Algebraic and
Geometric Topology at}
\immediate\write\gtoutfile{\s\s\s  http://www.maths.warwick.ac.uk/agt/AGTVol\thevolumenumber/agt-\thevolumenumber-\thepapernumber.abs.html}
\immediate\write\gtoutfile{\noexpand\\}
\immediate\write\gtoutfile{}
\ifx\theasciiabstract\relax
\immediate\write\gtoutfile{\theabstract}\else
\immediate\write\gtoutfile{\theasciiabstract}\fi
\immediate\write\gtoutfile{}
\immediate\write\gtoutfile{\noexpand\\}
\immediate\write\gtoutfile{}
\immediate\closeout\gtoutfile}}  

\def\maketitlepage{\makeagttitle\makeheadfile}
\let\makeshorttitle\maketitlepage
\let\maketitle\maketitlepage

\def\ifplaintex{\expandafter\ifx\csname documentclass\endcsname\relax}

\def\gtp{{\mathsurround=0pt\it $\cal G\mskip-2mu$eometry \&\ 
$\cal T\!\!$opology $\cal P\!$ublications}}  

\def\recd{{\small Received:\qua\receiveddate\ifx\reviseddate\relax
\else\qquad Revised:\qua\reviseddate\fi\par}} 

\def\lognumber#1{\def\thelognumber{#1}}
\def\volumenumber#1{\def\thevolumenumber{#1}}
\def\volumeyear#1{\def\thevolumeyear{#1}}
\def\papernumber#1{\def\thepapernumber{#1}}
\def\pagenumbers#1#2{\def\startpage{#1}\def\finishpage{#2}}
\def\published#1{\def\publishdate{#1}}

\def\received#1{\def\receiveddate{#1}}
\def\revised#1{\def\reviseddate{#1}}
\def\accepted#1{\def\accepteddate{#1}}

\def\asciiaddress#1{\def\theasciiaddress{#1}}

\let\\\par\let\thelognumber\relax\let\thevolumenumber\relax
\let\thepapernumber\relax\let\thevolumeyear\relax\let\startpage\relax
\let\finishpage\relax\let\publishdate\relax\let\receiveddate\relax
\let\reviseddate\relax\let\accepteddate\relax\let\theasciititle\relax
\let\theasciiauthors\relax\let\theasciiaddress\relax
\let\theasciiabstract\relax

\let\theasciiemail\relax

\ifplaintex
\font\logobig=cmssbx10 scaled 3836
\font\logomed=cmssbx10 scaled 2557
\else
\font\logobig=cmssbx10 scaled 4200
\font\logomed=cmssbx10 scaled 2800
\fi

\long\def\makeagttitle{   
\count0=\startpage
\agt\hfill      
\hbox to 45truept{\vbox to 0pt{\vglue -13truept{\logomed A\kern -.37em{\logobig 
T}\kern -.38em G}\vss}\hss}
\break
{\small Volume \thevolumenumber\ (\thevolumeyear)
\startpage--\finishpage\nl
Published: \publishdate}

\vglue .25truein

{\parskip=0pt\leftskip 0pt plus
1fil\def\\{\par\smallskip}{\Large\bf\thetitle}\par\medskip} \vglue
0.05truein

{\parskip=0pt\leftskip 0pt plus 1fil\def\\{\par}{\sc\theauthors}
\par\medskip}%
 
\vglue 0.03truein 

{\small\leftskip 25truept\rightskip 25truept{\bf Abstract}\stdspace\theabstract

{\bf AMS Classification}\stdspace\theprimaryclass
\ifx\thesecondaryclass\relax\else; \thesecondaryclass\fi\par
{\bf Keywords}\stdspace \thekeywords\par}\vglue 7truept

}   

\ifplaintex
\font\phead=cmsl9 scaled 950
\font\pnum=cmbx10 scaled 913
\font\pfoot=cmsl9 scaled 950
\headline{\vbox to 0pt{\vskip -4.5mm\line{\small\phead\ifnum
\count0=\startpage ISSN 1472-2739 (on-line) 1472-2747 (printed)
\hfill {\pnum\folio}\else\ifodd\count0\def\\{ }%
\ifx\theshorttitle\relax\thetitle\else\theshorttitle\fi\hfill{\pnum\folio}
\else\def\\{ and }{\pnum\folio}\hfill\ifx\theshortauthors\relax\theauthors
\else\theshortauthors\fi\fi\fi}\vss}}
\footline{\vbox to 0pt{\vglue 0mm\line{\small\pfoot\ifnum\count0=\startpage
\copyright\ \gtp\hfill\else
\agt, Volume \thevolumenumber\ (\thevolumeyear)\hfill\fi}\vss}}
\else
\font=cmsl9 scaled 1050
\font\lnum=cmbx10 
\font=cmsl9 scaled 1050
\makeatletter
\def\@oddhead{{\small\lhead\ifnum\count0=\startpage ISSN 1472-2739 
(on-line) 1472-2747 (printed)\hfill {\lnum\number\count0}\else\ifodd\count0
\def\\{ }\ifx\theshorttitle\relax \thetitle \else\theshorttitle\fi\hfill
{\lnum\number\count0}\else\def\\{ and }{\lnum\number\count0}
\hfill\ifx\theshortauthors\relax 
\theauthors\else\theshortauthors\fi\fi\fi}}\def\@evenhead{\@oddhead}
\def\@oddfoot{\small\lfoot\ifnum\count0=\startpage\copyright\ \gtp\hfill\else
\agt, Volume \thevolumenumber\ (\thevolumeyear)\hfill\fi}
\def\@evenfoot{\@oddfoot}
\makeatother
\fi
\let\maketitlepage\makeagttitle
\let\makeshorttitle\maketitlepage
\let\maketitle\maketitlepage

\newwrite\gtoutfile
\long\gdef\makeheadfile{  
{\def\\{, }\def\s{ }
\immediate\openout\gtoutfile head.xxx
\immediate\write\gtoutfile{To: math@arxiv.org}
\immediate\write\gtoutfile{Subject: put OR rep NNNNN:ppppp}
\immediate\write\gtoutfile{--text follows this line--}
\immediate\write\gtoutfile{Proxy-for: \ifx\theasciiauthors\relax
\theauthors\else\theasciiauthors\fi\s<\ifx\theasciiemail\relax\theemail\else\theasciiemail\fi>}
\immediate\write\gtoutfile{\noexpand\\}
\immediate\write\gtoutfile{Authors: \ifx\theasciiauthors\relax
\theauthors\else\theasciiauthors\fi}
{\def\\{ }\immediate\write\gtoutfile{Title: \ifx\theasciititle\relax
\thetitle\else\theasciititle\fi}}
\immediate\write\gtoutfile{Subj-class: GT or SG, GR etc}
\immediate\write\gtoutfile{MSC-class: \theprimaryclass\ifx\thesecondaryclass\relax\else, \thesecondaryclass\fi}
\immediate\write\gtoutfile{Journal-ref: Algebr. Geom. Topol. \thevolumenumber\s
(\thevolumeyear) \startpage-\finishpage}
\immediate\write\gtoutfile{Comments: Published by Algebraic and
Geometric Topology at}
\immediate\write\gtoutfile{\s\s\s  http://www.maths.warwick.ac.uk/agt/AGTVol\thevolumenumber/agt-\thevolumenumber-\thepapernumber.abs.html}
\immediate\write\gtoutfile{\noexpand\\}
\immediate\write\gtoutfile{}
\ifx\theasciiabstract\relax
\immediate\write\gtoutfile{\theabstract}\else
\immediate\write\gtoutfile{\theasciiabstract}\fi
\immediate\write\gtoutfile{}
\immediate\write\gtoutfile{\noexpand\\}
\immediate\write\gtoutfile{}
\immediate\closeout\gtoutfile}}  

\def\maketitlepage{\makeagttitle\makeheadfile}
\let\makeshorttitle\maketitlepage
\let\maketitle\maketitlepage

\def\ifplaintex{\expandafter\ifx\csname documentclass\endcsname\relax}

\def\gtp{{\mathsurround=0pt\it $\cal G\mskip-2mu$eometry \&\ 
$\cal T\!\!$opology $\cal P\!$ublications}}  

\def\recd{{\small Received:\qua\receiveddate\ifx\reviseddate\relax
\else\qquad Revised:\qua\reviseddate\fi\par}} 

\def\lognumber#1{\def\thelognumber{#1}}
\def\volumenumber#1{\def\thevolumenumber{#1}}
\def\volumeyear#1{\def\thevolumeyear{#1}}
\def\papernumber#1{\def\thepapernumber{#1}}
\def\pagenumbers#1#2{\def\startpage{#1}\def\finishpage{#2}}
\def\published#1{\def\publishdate{#1}}

\def\received#1{\def\receiveddate{#1}}
\def\revised#1{\def\reviseddate{#1}}
\def\accepted#1{\def\accepteddate{#1}}

\def\asciiaddress#1{\def\theasciiaddress{#1}}

\let\\\par\let\thelognumber\relax\let\thevolumenumber\relax
\let\thepapernumber\relax\let\thevolumeyear\relax\let\startpage\relax
\let\finishpage\relax\let\publishdate\relax\let\receiveddate\relax
\let\reviseddate\relax\let\accepteddate\relax\let\theasciititle\relax
\let\theasciiauthors\relax\let\theasciiaddress\relax
\let\theasciiabstract\relax

\let\theasciiemail\relax

\ifplaintex
\font\logobig=cmssbx10 scaled 3836
\font\logomed=cmssbx10 scaled 2557
\else
\font\logobig=cmssbx10 scaled 4200
\font\logomed=cmssbx10 scaled 2800
\fi

\long\def\makeagttitle{   
\count0=\startpage
\agt\hfill      
\hbox to 45truept{\vbox to 0pt{\vglue -13truept{\logomed A\kern -.37em{\logobig 
T}\kern -.38em G}\vss}\hss}
\break
{\small Volume \thevolumenumber\ (\thevolumeyear)
\startpage--\finishpage\nl
Published: \publishdate}

\vglue .25truein

{\parskip=0pt\leftskip 0pt plus
1fil\def\\{\par\smallskip}{\Large\bf\thetitle}\par\medskip} \vglue
0.05truein

{\parskip=0pt\leftskip 0pt plus 1fil\def\\{\par}{\sc\theauthors}
\par\medskip}%
 
\vglue 0.03truein 

{\small\leftskip 25truept\rightskip 25truept{\bf Abstract}\stdspace\theabstract

{\bf AMS Classification}\stdspace\theprimaryclass
\ifx\thesecondaryclass\relax\else; \thesecondaryclass\fi\par
{\bf Keywords}\stdspace \thekeywords\par}\vglue 7truept

}   

\ifplaintex
\font\phead=cmsl9 scaled 950
\font\pnum=cmbx10 scaled 913
\font\pfoot=cmsl9 scaled 950
\headline{\vbox to 0pt{\vskip -4.5mm\line{\small\phead\ifnum
\count0=\startpage ISSN 1472-2739 (on-line) 1472-2747 (printed)
\hfill {\pnum\folio}\else\ifodd\count0\def\\{ }%
\ifx\theshorttitle\relax\thetitle\else\theshorttitle\fi\hfill{\pnum\folio}
\else\def\\{ and }{\pnum\folio}\hfill\ifx\theshortauthors\relax\theauthors
\else\theshortauthors\fi\fi\fi}\vss}}
\footline{\vbox to 0pt{\vglue 0mm\line{\small\pfoot\ifnum\count0=\startpage
\copyright\ \gtp\hfill\else
\agt, Volume \thevolumenumber\ (\thevolumeyear)\hfill\fi}\vss}}
\else
\font=cmsl9 scaled 1050
\font\lnum=cmbx10 
\font=cmsl9 scaled 1050
\makeatletter
\def\@oddhead{{\small\lhead\ifnum\count0=\startpage ISSN 1472-2739 
(on-line) 1472-2747 (printed)\hfill {\lnum\number\count0}\else\ifodd\count0
\def\\{ }\ifx\theshorttitle\relax \thetitle \else\theshorttitle\fi\hfill
{\lnum\number\count0}\else\def\\{ and }{\lnum\number\count0}
\hfill\ifx\theshortauthors\relax 
\theauthors\else\theshortauthors\fi\fi\fi}}\def\@evenhead{\@oddhead}
\def\@oddfoot{\small\lfoot\ifnum\count0=\startpage\copyright\ \gtp\hfill\else
\agt, Volume \thevolumenumber\ (\thevolumeyear)\hfill\fi}
\def\@evenfoot{\@oddfoot}
\makeatother
\fi
\let\maketitlepage\makeagttitle
\let\makeshorttitle\maketitlepage
\let\maketitle\maketitlepage

\newwrite\gtoutfile
\long\gdef\makeheadfile{  
{\def\\{, }\def\s{ }
\immediate\openout\gtoutfile head.xxx
\immediate\write\gtoutfile{To: math@arxiv.org}
\immediate\write\gtoutfile{Subject: put OR rep NNNNN:ppppp}
\immediate\write\gtoutfile{--text follows this line--}
\immediate\write\gtoutfile{Proxy-for: \ifx\theasciiauthors\relax
\theauthors\else\theasciiauthors\fi\s<\ifx\theasciiemail\relax\theemail\else\theasciiemail\fi>}
\immediate\write\gtoutfile{\noexpand\\}
\immediate\write\gtoutfile{Authors: \ifx\theasciiauthors\relax
\theauthors\else\theasciiauthors\fi}
{\def\\{ }\immediate\write\gtoutfile{Title: \ifx\theasciititle\relax
\thetitle\else\theasciititle\fi}}
\immediate\write\gtoutfile{Subj-class: GT or SG, GR etc}
\immediate\write\gtoutfile{MSC-class: \theprimaryclass\ifx\thesecondaryclass\relax\else, \thesecondaryclass\fi}
\immediate\write\gtoutfile{Journal-ref: Algebr. Geom. Topol. \thevolumenumber\s
(\thevolumeyear) \startpage-\finishpage}
\immediate\write\gtoutfile{Comments: Published by Algebraic and
Geometric Topology at}
\immediate\write\gtoutfile{\s\s\s  http://www.maths.warwick.ac.uk/agt/AGTVol\thevolumenumber/agt-\thevolumenumber-\thepapernumber.abs.html}
\immediate\write\gtoutfile{\noexpand\\}
\immediate\write\gtoutfile{}
\ifx\theasciiabstract\relax
\immediate\write\gtoutfile{\theabstract}\else
\immediate\write\gtoutfile{\theasciiabstract}\fi
\immediate\write\gtoutfile{}
\immediate\write\gtoutfile{\noexpand\\}
\immediate\write\gtoutfile{}
\immediate\closeout\gtoutfile}}  

\def\maketitlepage{\makeagttitle\makeheadfile}
\let\makeshorttitle\maketitlepage
\let\maketitle\maketitlepage

\def\ifplaintex{\expandafter\ifx\csname documentclass\endcsname\relax}

\def\gtp{{\mathsurround=0pt\it $\cal G\mskip-2mu$eometry \&\ 
$\cal T\!\!$opology $\cal P\!$ublications}}  

\def\recd{{\small Received:\qua\receiveddate\ifx\reviseddate\relax
\else\qquad Revised:\qua\reviseddate\fi\par}} 

\def\lognumber#1{\def\thelognumber{#1}}
\def\volumenumber#1{\def\thevolumenumber{#1}}
\def\volumeyear#1{\def\thevolumeyear{#1}}
\def\papernumber#1{\def\thepapernumber{#1}}
\def\pagenumbers#1#2{\def\startpage{#1}\def\finishpage{#2}}
\def\published#1{\def\publishdate{#1}}

\def\received#1{\def\receiveddate{#1}}
\def\revised#1{\def\reviseddate{#1}}
\def\accepted#1{\def\accepteddate{#1}}

\def\asciiaddress#1{\def\theasciiaddress{#1}}

\let\\\par\let\thelognumber\relax\let\thevolumenumber\relax
\let\thepapernumber\relax\let\thevolumeyear\relax\let\startpage\relax
\let\finishpage\relax\let\publishdate\relax\let\receiveddate\relax
\let\reviseddate\relax\let\accepteddate\relax\let\theasciititle\relax
\let\theasciiauthors\relax\let\theasciiaddress\relax
\let\theasciiabstract\relax

\let\theasciiemail\relax

\ifplaintex
\font\logobig=cmssbx10 scaled 3836
\font\logomed=cmssbx10 scaled 2557
\else
\font\logobig=cmssbx10 scaled 4200
\font\logomed=cmssbx10 scaled 2800
\fi

\long\def\makeagttitle{   
\count0=\startpage
\agt\hfill      
\hbox to 45truept{\vbox to 0pt{\vglue -13truept{\logomed A\kern -.37em{\logobig 
T}\kern -.38em G}\vss}\hss}
\break
{\small Volume \thevolumenumber\ (\thevolumeyear)
\startpage--\finishpage\nl
Published: \publishdate}

\vglue .25truein

{\parskip=0pt\leftskip 0pt plus
1fil\def\\{\par\smallskip}{\Large\bf\thetitle}\par\medskip} \vglue
0.05truein

{\parskip=0pt\leftskip 0pt plus 1fil\def\\{\par}{\sc\theauthors}
\par\medskip}%
 
\vglue 0.03truein 

{\small\leftskip 25truept\rightskip 25truept{\bf Abstract}\stdspace\theabstract

{\bf AMS Classification}\stdspace\theprimaryclass
\ifx\thesecondaryclass\relax\else; \thesecondaryclass\fi\par
{\bf Keywords}\stdspace \thekeywords\par}\vglue 7truept

}   

\ifplaintex
\font\phead=cmsl9 scaled 950
\font\pnum=cmbx10 scaled 913
\font\pfoot=cmsl9 scaled 950
\headline{\vbox to 0pt{\vskip -4.5mm\line{\small\phead\ifnum
\count0=\startpage ISSN 1472-2739 (on-line) 1472-2747 (printed)
\hfill {\pnum\folio}\else\ifodd\count0\def\\{ }%
\ifx\theshorttitle\relax\thetitle\else\theshorttitle\fi\hfill{\pnum\folio}
\else\def\\{ and }{\pnum\folio}\hfill\ifx\theshortauthors\relax\theauthors
\else\theshortauthors\fi\fi\fi}\vss}}
\footline{\vbox to 0pt{\vglue 0mm\line{\small\pfoot\ifnum\count0=\startpage
\copyright\ \gtp\hfill\else
\agt, Volume \thevolumenumber\ (\thevolumeyear)\hfill\fi}\vss}}
\else
\font=cmsl9 scaled 1050
\font\lnum=cmbx10 
\font=cmsl9 scaled 1050
\makeatletter
\def\@oddhead{{\small\lhead\ifnum\count0=\startpage ISSN 1472-2739 
(on-line) 1472-2747 (printed)\hfill {\lnum\number\count0}\else\ifodd\count0
\def\\{ }\ifx\theshorttitle\relax \thetitle \else\theshorttitle\fi\hfill
{\lnum\number\count0}\else\def\\{ and }{\lnum\number\count0}
\hfill\ifx\theshortauthors\relax 
\theauthors\else\theshortauthors\fi\fi\fi}}\def\@evenhead{\@oddhead}
\def\@oddfoot{\small\lfoot\ifnum\count0=\startpage\copyright\ \gtp\hfill\else
\agt, Volume \thevolumenumber\ (\thevolumeyear)\hfill\fi}
\def\@evenfoot{\@oddfoot}
\makeatother
\fi
\let\maketitlepage\makeagttitle
\let\makeshorttitle\maketitlepage
\let\maketitle\maketitlepage

\newwrite\gtoutfile
\long\gdef\makeheadfile{  
{\def\\{, }\def\s{ }
\immediate\openout\gtoutfile head.xxx
\immediate\write\gtoutfile{To: math@arxiv.org}
\immediate\write\gtoutfile{Subject: put OR rep NNNNN:ppppp}
\immediate\write\gtoutfile{--text follows this line--}
\immediate\write\gtoutfile{Proxy-for: \ifx\theasciiauthors\relax
\theauthors\else\theasciiauthors\fi\s<\ifx\theasciiemail\relax\theemail\else\theasciiemail\fi>}
\immediate\write\gtoutfile{\noexpand\\}
\immediate\write\gtoutfile{Authors: \ifx\theasciiauthors\relax
\theauthors\else\theasciiauthors\fi}
{\def\\{ }\immediate\write\gtoutfile{Title: \ifx\theasciititle\relax
\thetitle\else\theasciititle\fi}}
\immediate\write\gtoutfile{Subj-class: GT or SG, GR etc}
\immediate\write\gtoutfile{MSC-class: \theprimaryclass\ifx\thesecondaryclass\relax\else, \thesecondaryclass\fi}
\immediate\write\gtoutfile{Journal-ref: Algebr. Geom. Topol. \thevolumenumber\s
(\thevolumeyear) \startpage-\finishpage}
\immediate\write\gtoutfile{Comments: Published by Algebraic and
Geometric Topology at}
\immediate\write\gtoutfile{\s\s\s  http://www.maths.warwick.ac.uk/agt/AGTVol\thevolumenumber/agt-\thevolumenumber-\thepapernumber.abs.html}
\immediate\write\gtoutfile{\noexpand\\}
\immediate\write\gtoutfile{}
\ifx\theasciiabstract\relax
\immediate\write\gtoutfile{\theabstract}\else
\immediate\write\gtoutfile{\theasciiabstract}\fi
\immediate\write\gtoutfile{}
\immediate\write\gtoutfile{\noexpand\\}
\immediate\write\gtoutfile{}
\immediate\closeout\gtoutfile}}  

\def\maketitlepage{\makeagttitle\makeheadfile}
\let\makeshorttitle\maketitlepage
\let\maketitle\maketitlepage

\lognumber{4}
\volumenumber{2}
\volumeyear{2002}
\papernumber{4}
\published{1 February 2002}
\pagenumbers{37}{50}
\received{19 November 2001}
\revised{16 January 2002}
\accepted{27 January 2002}

\usepackage{amssymb} 
\usepackage{amsmath,amscd} 
\usepackage{epsf}
\usepackage{latexsym}
\def\S{section }

\newtheorem{theorem}{Theorem}[section]
\newtheorem{lemma}[theorem]{Lemma}
\newtheorem{corollary}[theorem]{Corollary}
\newtheorem{conjecture}[theorem]{Conjecture}
\newtheorem{scholium}[theorem]{Scholium}

\theoremstyle{remark}
\newtheorem{example}[theorem]{Example}
\newtheorem*{rem}{Remark}

\def\Aut{\mbox{\rm{Aut}}} 
\def\Mod{\mbox{\rm{Mod}}}  
 \def\fix{\mbox{\rm{fix}}}

\def\rk{\mbox{\rm{rk}}}   

\begin{document}
\title{The co-rank conjecture for 3-manifold groups}                    
\authors{Christopher J. Leininger\\Alan W. Reid}                  
\address{Department of Mathematics\\ University of Texas at Austin\\Austin, 
TX  78712-1082, USA\\\smallskip\\{\rm FAX: (512)-471-9038}}                  
\asciiaddress{Department of Mathematics\\University of Texas at 
Austin\\Austin, TX  78712-1082, USA\\FAX: (512)-471-9038}                  
\email{clein@math.utexas.edu, areid@math.utexas.edu }                     

\begin{abstract}  
In this paper we construct explicit examples of both closed and
non-compact finite volume hyperbolic manifolds which provide
counterexamples to the conjecture that the co-rank of a 3-manifold
group (also known as the cut number) is bounded below by one-third the
first Betti number.
\end{abstract}


\maketitle

\section{Introduction}

The {\em co-rank} of a group $G$, which we denote by $c_{1}(G)$, is
the maximal rank of a free group homomorphically surjected by $G$.
Clearly, $c_{1}(G) \leq b_{1}(G) = \dim(G^{ab}\otimes {\mathbb Q})$,
where $G^{ab}\! =\! (G/[G,G])$. If $M$ is a manifold, $c_{1}(M) =
c_{1}(\pi_{1}(M))$.
If $M$ is a compact 3-manifold, $c_1(M)$ is
also called the {\em cut number} of $M$, and is equal to the maximal
number of
components of a surface $F$ embedded in $M$ for which $M\setminus F$
is connected.

The free Abelian groups ${\mathbb Z}^{n}$ show that, in general, there
is no lower bound for $c_{1}(G)$ in terms of $b_{1}(G)$. For a genus
$g$ surface group it is well-known that $c_{1}(\Sigma_{g})=g$ (as is
proved below).  In his talk \cite{Mont}, J Stallings discussed the
following conjecture on a lower bound for $c_1(M)$ for $M$ a compact
3-manifold, which has recently received some attention. According to
A Sikora this conjecture has its origins in work of T Kerler
connected to quantum invariants of 3-manifolds.

\begin{conjecture}
\label{nochance}
If $M$ is a $3$-manifold, then $c_{1}(M) \geq \frac{b_{1}(M)}{3}$.
\end{conjecture}

Notice that as particular cases,
if $b_{1}(M) = 4$ or $5$, Conjecture \ref{nochance} would imply
$c_{1}(M) \geq 2$.
In this note we construct explicit counterexamples to this conjecture.
In particular, we prove:

\begin{theorem} \label{main}

{\rm(1)}\qua There exist closed hyperbolic $3$-manifolds $M$ such that $b_{1}(M)=5$ and $c_{1}(M)=1$.

{\rm(2)}\qua There exist compact $3$-manifolds $M$ with toroidal boundary so that
$b_{1}(M)$\break $=4$ and $c_{1}(M)=1$.
\end{theorem}

\begin{rem}A Sikora \cite{Sik} has also recently proved the existence
of counterexamples to this conjecture
with $b_{1}\geq 7$ and $c_{1} \leq 2$. Moreover, S. Harvey 
\cite{Ha} has recently constructed examples of closed hyperbolic 3-manifolds
with $b_1$ arbitrarily large, yet $c_1 = 1$.\end{rem}

The rest of the paper is organized as follows.  In \S 2 we fix
some notation and make a few elementary observations about co-rank and
surface groups.  In \S 3 we discuss automorphisms of surfaces and
when they extend over a handlebody. This is needed to prove Theorem
\ref{main}.
We discuss the contents of \S 3 further in \S 6. Parts (1) and (2) of
Theorem \ref{main} are proved in sections 4 and 5.

\medskip

{\bf Acknowledgement}\qua The authors wish to thank C\thinspace McA Gordon
and D\thinspace D Long for useful conversations regarding this work,
K Johannson for providing us a copy of his unpublished manuscript
\cite{JJ}, and particularly the organizers of the Workshop in Groups
and 3-Manifolds, CRM Montreal, June/July 2001 for support, where they
first became aware of this problem.

\noindent
This work was partially supported by the NSF.

\section{Background and notation: Co-rank and surfaces}

\subsection{Notation}

Throughout the rest of this paper $F_{k}$ will denote a free group of rank
$k$, $\Sigma_{g}$ will denote a closed, oriented surface of genus $g$, and
$\Mod(\Sigma_{g})$ will denote its mapping class group.
That is, $\Mod(\Sigma_{g})$ is the group of isotopy classes of orientation
preserving diffeomorphisms of $\Sigma_{g}$.
We will refer to an element of $\Mod(\Sigma_{g})$, or any representative
of
that element, as an {\em automorphism} of $\Sigma_{g}$.

The canonical homomorphism $\Mod(\Sigma_{g}) \rightarrow
\Aut(H_{1}(\Sigma_{g}))$ has a non-trivial kernel which is called the {\em
Torelli group}.
We will denote it by 
$$ {\cal I}(\Sigma_{g}) \triangleleft \Mod(\Sigma_{g}).$$
An automorphism $f$ is in ${\cal I}(\Sigma_{g})$ if and only if
$f_{*}:H_{1}(\Sigma_{g}) \rightarrow H_{1}(\Sigma_{g})$ is the identity.

\subsection{Surfaces}

In this section we record some lemmas concerning epimorphisms from
surface groups to free groups. Since we shall make use
of it, we give the proof that $c_1(\Sigma_g) = g$, although this is 
well-known. 

\begin{lemma} \label{c1surface}
$c_{1}(\Sigma_{g})=g$.
\end{lemma}

\noindent
\proof Let $\phi: \pi_{1}(\Sigma_{g}) \rightarrow F_{k}$ be a surjective homomorphism and let $f: \Sigma_{g} \rightarrow X_{k}$ be a map where $X_{k}$ is a wedge of $k$ circles with $\pi_{1}(X_{k})$ identified with $F_{k}$ so that $f_{*} = \phi$.
Making $f$ transverse to $k$ points (one in each circle, each different than the wedge point), the preimage is a disjoint union of $k$ 1-submanifolds $\gamma_{1},...,\gamma_{k}$.
It is easy to see that $[\gamma_{1}],...,[\gamma_{k}]$ must represent $k$ linearly independent elements of $H_{1}(\Sigma_{g};{\mathbb Z})$ (with a pull-back orientation).
Since the curves $\gamma_{1},...,\gamma_{k}$ are pairwise disjoint, the intersection form on $H_{1}(\Sigma_{g};{\mathbb Z})$ is trivial on the span $<[\gamma_{1}],...,[\gamma_{k}]>$.
This is a non-degenerate, skew-symmetric form, so that $k \leq \frac{1}{2}b_{1}(\Sigma_{g}) = g$.

Clearly $c_{1}(\Sigma_{g}) \geq g$, so $c_{1}(\Sigma_{g}) = g$.
\endproof

We will also need the following well-known fact, whose proof we sketch.

\begin{lemma} \label{surfacehandlebody}
Suppose $\phi:\pi_{1}(\Sigma_{g}) \rightarrow F_{g}$ is an epimorphism.
Then there exists a handlebody $H$ with $\pi_{1}(H)$ identified to $F_{g}$, and a homeomorphism $\tau: \Sigma_{g} \rightarrow \partial H$, such that if $i:\partial H \rightarrow H$ is the inclusion map, then $(i \circ \tau)_{*} = \phi$.
\end{lemma}

\proof Given $\phi: \pi_{1}(\Sigma_{g}) \rightarrow F_{g}$, let $f:\Sigma_{g} \rightarrow X_{g}$ and $\gamma_{1},...,\gamma_{g}$ be as in the previous proof.
We may homotope $f$ so that no component of $\gamma_{j}$ is
homotopically trivial for any $j=1,...,g$.
We construct a compression body $H$ by first thickening $\Sigma_{g}$ to
$\Sigma_{g} \times [-1,1]$, then attaching $2$-handles along the curves on
$\Sigma_{g} \times \{-1 \}$ corresponding to $\cup_{j=1}^{g}\gamma_{j}$ on $\Sigma_{g}$,
and finally capping off all $2-$sphere boundary components with
$3-$handles.

It must be that $H$ is a handlebody.
If this were not the case, one could find an (oriented) curve $\gamma$ on $\Sigma_{g}$, disjoint from $\gamma_{1},...,\gamma_{g}$ with $[\gamma],[\gamma_{1}],...,[\gamma_{g}]$ being linearly independent.
As in the previous proof, this is impossible by dimension considerations.

Let $\tau: \Sigma_{g} \rightarrow \Sigma_{g} \times \{ 1 \} = \partial H$ be the obvious homeomorphism.
By construction, we can extend $f$ to $\widetilde{f}:H \rightarrow X_{g}$ (that is, $\widetilde{f} \circ \tau = f$).
Moreover, the induced homomorphism $\widetilde{f}_{*}:\pi_{1}(H) \rightarrow \pi_{1}(X_{g}) = F_{g}$ must then be surjective.
It follows from the Hopfian property for free groups (see \cite{MKS}) that $\widetilde{f}_{*}$ must be an isomorphism.
Therefore $\widetilde{f}_{*}$ identifies $\pi_{1}(H)$ with $F_{k}$ and $(i \circ \tau)_{*} = f_{*} = \phi$, where $i:\partial H \rightarrow H$ is the inclusion.\endproof

\section{Extending automorphisms of surfaces}

Recall that an automorphism $f:\Sigma_{g} \rightarrow \Sigma_{g}$ extends
over a handlebody $H$ if there exists a homeomorphism $\tau:\Sigma_{g} \rightarrow \partial H$ and a homeomorphism $\tilde{f}:H \rightarrow H$ such that $\tau^{-1} \circ \tilde{f}|_{\partial H} \circ \tau = f$. 

\begin{lemma} \label{extendiff}
For $f \in \Mod(\Sigma_{g})$, $f$ extends over a handlebody $H$ if and
only if for every simple closed curve $\gamma \subset \Sigma_{g}$ for which $\tau(\gamma)$ bounds
a disk in $H$, $\tau(f(\gamma))$ also bounds a disk in $H$.
\end{lemma}

\proof If $f$ extends over a handlebody, then clearly $\tau(\gamma)$
bounds a disk if and only if $\tau(f(\gamma))$ does.

Suppose that for every simple closed curve $\gamma \subset \Sigma_{g}$
for which $\tau(\gamma)$ bounds a disk in $H$, $\tau(f(\gamma))$ also bounds a disk.
Choose a complete set of (pairwise disjoint) meridional disks
$D_{1},...,D_{g}$ for $H$.
The images $\tau \circ f \circ \tau^{-1} (\partial(D_{1})),...,\tau \circ f \circ \tau^{-1} (\partial(D_{g}))$ bound disks $D'_{1},...,D'_{g}$ which we may assume are pairwise disjoint.

We now extend $f$ to $\tilde{f}:H \rightarrow H$.
This is done by first extending over the disks $D_{j}$, mapping these to
$D'_{j}$ by any homeomorphism which extends $(\tau \circ f \circ \tau^{-1})|_{\partial D_{j}}$, and
then extending over a regular neighborhood of $\partial H \cup
\bigcup_{j=1}^{g} D_{j}$.
What is left is a $3$-ball, and we may extend over this in any fashion.
\endproof

The important point for us is existence of certain types of automorphisms
that do not extend.
The following theorem was proven independently by Johannson and Johnson
\cite{JJ} and Casson \cite{Ca}.  
Neither of these works were ever published, so we include in \S 6 a sketch
of the proof (as given in \cite{JJ}) for completeness.

\begin{theorem} \label{noextend}
For every $g \geq 2$, there exists $f \in {\cal I}(\Sigma_{g})$ which does
not extend over any handlebody. Furthermore, for every odd integer $n$,
$f^n$ does not extend.
\end{theorem}

Indeed, using Thurston's classification of automorphisms of surfaces
\cite{Thaut}, the automorphism $f$ can be chosen to be pseudo-Anosov. We
state this below and prove this in \S 6, since we make use of some of the
notation developed in proving Theorem \ref{noextend}.

\begin{theorem} \label{pA}
For any $g \geq 2$, there exists pseudo-Anosov mapping classes $f \in
{\cal I}(\Sigma_{g})$ so that no odd power of $f$ extends over any
handlebody.
\end{theorem}

\rem In the genus $2$ case, {\em every} automorphism
$f \in {\cal I}(\Sigma_{2})$ which does not extend over a handlebody is
pseudo-Anosov.
To see this, first note that by Thurston's classification of automorphisms of surfaces, and because ${\cal I}(\Sigma_{2})$ is torsion free, $f$ is either reducible or pseudo-Anosov.
If $f \in {\cal I }(\Sigma_{2})$ is reducible, we can extend $f$ to $\widetilde{f} : H \rightarrow H$ where $H$ is a compression body defined by the reducing curves of $f$ with upper boundary, $\partial_{+}H$, identified to $\Sigma_{2}$.
The lower boundary of the compression body, $\partial_{-}H$, must be a
(possibly empty) disjoint union of tori.
Since the map induced by inclusion $i_{*}:H_{1}(\partial_{-}H) \rightarrow H_{1}(H)$ is injective and $\widetilde{f}$ acts trivially on $H_{1}(H)$, it follows that $(\widetilde{f}|_{\partial_{-}H})_{*}:H_{1}(\partial_{-}H) \rightarrow H_{1}(\partial_{-}H)$ must be the identity.
An automorphism acting trivially on the homology of a torus is trivial, and so we can compress away the lower boundary completely and extend.
It follows that if $f \in {\cal I}(\Sigma_{2})$ does not extend over any handlebody, then it must be pseudo-Anosov.

\section{Examples: Closed 3-manifolds}

Let $M_{f}$ be the mapping torus of a pseudo-Anosov automorphism $f \in
{\cal I}(\Sigma_{2})$ that does not extend over any handlebody (by Theorem
\ref{pA}, such an $f$ exists).

\begin{theorem} \label{mainclosed}
$M_{f}$ is hyperbolic, $b_{1}(M_{f}) = 5$, and $c_{1}(M_{f})=1$
\end{theorem}

\proof Since $f$ is pseudo-Anosov, $M_{f}$ is hyperbolic by
Thurston's Geometerization Theorem for Haken Manifolds \cite{Tfiber}.
An elementary calculation shows that $b_{1}(M_{f}) = 1+ \rk(\fix(f_{*}))$
where $\fix(f_{*})$ is the fixed subgroup of $H_{1}(\Sigma_{2};{\mathbb
Z})$.
So, $b_{1}(M_{f}) = 5$.

Suppose $c_{1}(M_{f}) > 1$.
Then there exists an epimorphism $\phi:\pi_{1}(M_{f}) \rightarrow F_{2}$.
Let $\Sigma_{2} \subset M_{f}$ denote the fiber, and note that
$\pi_{1}(\Sigma_{2}) \triangleleft \pi_{1}(M_{f})$.
This implies that $\phi(\pi_{1}(\Sigma_{2})) \triangleleft F_{2}$ is a
finitely generated normal subgroup, which must therefore be either trivial
or of finite index.
Now $\phi(\pi_{1}(\Sigma_{2}))$ cannot be trivial since this would imply
$F_{2}$ is the image of $\pi_{1}(M_{f}) / \pi_{1}(\Sigma_{2}) \cong
{\mathbb Z}$.
It follows that $\phi(\pi_{1}(\Sigma_{2}))$ must have finite index.
Note that $rk(\phi(\pi_{1}(\Sigma_{2}))) = 1+ [F_{2}:\phi(\pi_{1}(\Sigma_{2}))]$.
By Lemma \ref{c1surface}, $c_{1}(\Sigma_{2}) = 2$, so that $[F_{2}:\phi(\pi_{1}(\Sigma_{2}))]=1$ and
$\phi|_{\pi_{1}(\Sigma_{2})}$ must be surjective.

By Lemma \ref{surfacehandlebody}, there exists a handlebody $H$ with
$\Sigma_{2} = \partial H$ such that the inclusion $i : \Sigma_{2}
\rightarrow H$ has $i_{*} = \phi|_{\pi_{1}(\Sigma_{2})}$.
Since $f$ does not extend over any handlebody, Lemma \ref{extendiff}
implies the existence of a simple closed curve $\gamma \subset \Sigma_{2}$
such that $\gamma$ bounds a disk in $H$, but $f(\gamma)$ does not.
By Dehn's Lemma and the Loop Theorem (see eg \cite{Ro}), $f(\gamma)$ is
homotopically non-trivial in $H$.
Representing these curves by based loops of the same name in $\Sigma_{2}$
(with basepoint fixed by $f$), this says $i_{*}([\gamma]) = 1$, while
$i_{*}(f_{*}([\gamma])) = i_{*}([f(\gamma)]) \neq 1$.
That is, $\phi([\gamma])=1$ while $\phi(f_{*}([\gamma])) \neq 1$.

Now we note that $\pi_{1}(M_{f})$ is an HNN extension of $\pi_{1}(\Sigma_{2})$ with conjugation by the stable letter, $t$, acting by $f_{*}$.
It follows that
$$1 = \phi(t) \phi(t^{-1}) = \phi(t) \phi([\gamma]) \phi(t^{-1}) = \phi(t [\gamma] t^{-1}) = \phi(f_{*}[\gamma]) \neq 1.$$
This contradiction proves that $c_{1}(M_{f})=1$.\endproof

\begin{corollary}
\label{blob}
There are infinitely many closed hyperbolic $3$-manifolds $M$ with $b_{1}(M)=5$
and $c_{1}(M) =1$.
\end{corollary}

\proof Given $M_{f}$ as above, with $f$ as in Theorem \ref{pA}, the
cyclic covers $M_{f^{n}}$, for odd $n$, provide infinitely many such
manifolds.
\endproof

Using the nature of the construction of the automorphism $f$ in Theorem
\ref{mainclosed} we can extend Corollary \ref{blob} to the following; a
proof is given in \S 6.
\begin{corollary}
\label{newblob}
There are infinitely many non-commensurable 
closed hyperbolic $3$-manifolds $M$ with $b_{1}(M)=5$
and $c_{1}(M) =1$.
\end{corollary}
A similar argument can be made to work for genus $3$ bundles.
In this case if we consider an $f \in {\cal I}(\Sigma_3)$ as in Theorem
\ref{pA}, then $b_1(M_f) = 7$, so that
Conjecture \ref{nochance}
would predict $c_1(M_f) = 3$.
\begin{theorem}
\label{genus3}
There are infinitely many closed hyperbolic $3$-manifolds $M$ with $b_{1}(M) = 7$ and $c_{1}(M) \leq 2$.\end{theorem}

\proof Let $f \in {\cal I}(\Sigma_3)$ be as in
Theorem \ref{pA}. The proof of Theorem \ref{mainclosed}
applies verbatim to the bundle $M_{f}$.
The only point to remark being that
if $\pi_1(M_f)$ surjects a free group of rank $3$, then the fiber
group must also surject.
\endproof

\rem Notice that the argument breaks down for genus $4$
bundles.
In this case, $b_1(M_f) = 9$ so that a counterexample to Conjecture
\ref{nochance} requires $c_1(M_f) \leq 2$.
The argument above only guarantees $c_{1}(M_{f}) \leq 3$.

\section{Examples: Bounded $3$-manifolds}

Here we sketch the proof of the second part of Theorem \ref{main}.

\begin{theorem} \label{mainfinite}
There exists compact $3$-manifolds $M$ with toroidal boundary so that
$b_{1}(M)=4$ and $c_{1}(M)=1$.
\end{theorem}

\begin{rem}In fact, the manifolds we construct below can
be shown to be hyperbolic using Thurston's Geometerization Theorem for
Haken manifolds (see \cite{O} for example), and a variant of
\cite{Hcurve} (see \cite{CG} for a proof that $M$ is irreducible).
\end{rem}

\proof[Sketch of proof]
Let $\alpha,\beta,\gamma,\delta,$ and $\epsilon$ be the simple closed
curves shown on the surface $\Sigma_{2}$, and let $a,b,c,$ and $d$ be the
generators for $\pi_{1}(\Sigma_{2})$ shown (see Figure 1).

\begin{figure}[h]
\centerline{\epsfysize=1.1truein\epsfxsize=5.2truein\epsfbox{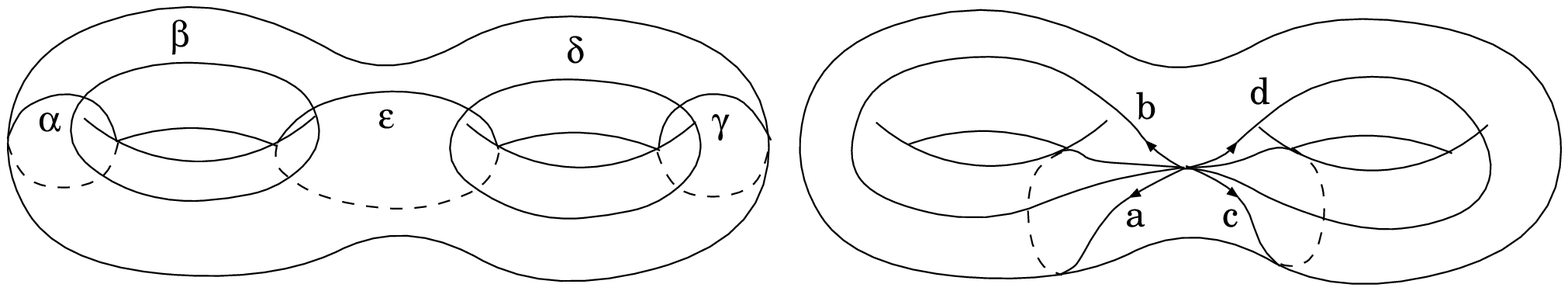}}
\nocolon\caption{}
\end{figure}

We denote the Dehn twist in $\alpha,\beta,\gamma,\delta,$ and $\epsilon$
by $T_{\alpha},T_{\beta},T_{\gamma},T_{\delta},$ and $T_{\epsilon}$
respectively.
Let $C_{1}$ denote the separating simple closed curve on $\Sigma_{2}$
represented by $[a,b] = aba^{-1}b^{-1}$, and set
$$C_{2} = T_{\alpha}^{-1} \circ (T_{\alpha} \circ T_{\epsilon} \circ
T_{\gamma}^{2} \circ T_{\beta}^{-1} \circ T_{\delta}^{-1})^{3}(C_{1})$$
Now construct a $3$-manifold $M$ by attaching $2$-handles to $\Sigma_{2}
\times [-1,1]$ along $C_{1}$ in $\Sigma_{2} \times \{ -1 \}$ and along
$C_{2}$ in $\Sigma_{2} \times \{ 1 \}$.
Since both $C_{1}$ and $C_{2}$ are separating curves, $b_{1}(M) =
b_{1}(\Sigma_{2}) = 4$.
We claim that $c_{1}(M) = 1$.

We begin by finding a presentation for $\pi_{1}(M)$.
By considering the action of each of the Dehn twists above on
$\pi_{1}(\Sigma)$, we can explicitly write down a word $w=w(a,b,c,d)$
representing $C_{2}$.
It is given by
\begin{gather*}
w=a c^{-1} a b c^{-1} a d c^{-3} a c^{-1} a b c^{-1} a b a c^{-1} a b
c^{-1} a d c^{-3} a c^{-1} a b\\
c^{-1} a d c^{-3} a d c^{-4} a d c^{-4} a d c^{-3} a c^{-1} a b c^{-1} a
d c^{-3} a c^{-1} a b a^{-1} b^{-1} a^{-1} c b^{-1} a^{-1} c a^{-1}\\
c^{3} d^{-1} a^{-1} c b^{-1} a^{-1} c a^{-1} b^{-1} a^{-1} c      b^{-1}
a^{-1} c a^{-1}  c^{3} d^{-1} a^{-1} c b^{-1} a^{-1} c a^{-1}\\
c^{3} d^{-1} a^{-1} c^{4} d^{-1} a^{-1} c^{4} d^{-1} a^{-1} c^{3} d^{-1}
a^{-1}    c b^{-1} a^{-1} c a^{-1} c^{3} d^{-1} a ^{-1} c.
\end{gather*}
An application of Van Kampen's Theorem implies
$$\pi_{1}(M) = <a,b,c,d \, | \, [a,b],[c,d],w>.$$
To prove that $c_{1}(M)=1$, we suppose there exists an epimorphism
$\phi: \pi_{1}(M) \rightarrow F_{2}$, and find a contradiction.
The idea of the proof is as follows.
Since $[a,b] = 1$ and $[c,d] = 1$ and
maximal Abelian subgroups of free groups are cyclic, there exists $g,h
\in F_{2}$ such that $\phi(a),\phi(b) \in <g>$ and $\phi(c),\phi(d) \in
<h>$.
We let $m,n,j,k \in {\mathbb Z}$ be such that
$$\phi(a) = g^{m} \, , \, \phi(b) = g^{n} \, , \, \phi(c) = h^{k} \, \,
\mbox{ and } \phi(d) = h^{j}.$$
The subgroup of $F_{2}$ generated by $g$ and $h$ contains the subgroup
generated by $\phi(a),\phi(b),\phi(c),$ and $\phi(d)$, but $\phi$ is
surjective, so $g$ and $h$ must generate all of $F_{2}$.
By the Hopfian property for free groups, we see that $g$ and $h$ form a
basis for $F_{2}$.
Since $\phi$ is a homomorphism, we have
$$1 = \phi(w(a,b,c,d)) = w(\phi(a),\phi(b),\phi(c),\phi(d)) =
w(g^{m},g^{n},h^{k},h^{j}).$$
This imposes restrictions on the integers $m,n,j,$ and $k$.
There are then several cases to analyze, each of which results in a
contradiction to the surjectivity of $\phi$.
It follows that there is no such epimorphism $\phi$, and hence
$c_{1}(M)=1$.
\endproof

\begin{rem} Another proof of this theorem goes as follows.
Choose two separating curves $C_{1}$ and $C_{2}$ on $\Sigma_{2}$ so that (a) $C_{1}
\cup C_{2}$ fills $\Sigma_{2}$, and (b) for every handlebody $H$ with
$\partial H = \Sigma_{2}$, at most one of $C_{1}$ and $C_{2}$ bounds a
disk.
If one constructs $M$ as in the above proof, then $M$ would provide an
example proving the theorem.
To see this, note that any epimorphism from $\pi_{1}(M)$ onto $F_{2}$
would induces a epimorphism from $\pi_{1}(\Sigma_{2})$ onto $F_{2}$ in
which $C_{1}$ and $C_{2}$ are both mapped to $1$.
By applying Dehn's Lemma and the Loop Theorem along with Lemma
\ref{surfacehandlebody} above, we would have a contradiction.\end{rem}

Of course, the difficulty is in finding two curves satisfying (a) and (b).
The two curves $C_{1}$ and $C_{2}$ in the given proof {\em do} satisfy
(a) and (b)--condition (b) is essentially what is shown in the proof
(which we only sketched).
In fact, Lemma 2.2 of \cite{CL} describes an algorithm to decide if two
curves can both bound disks in any handlebody, so it should possible to
implement this to give another proof that $C_{1}$ and $C_{2}$ satisfy (b).
This algorithm is based on analyzing the intersections of the pair of
curves, and it seems that likely that this could be computationally more
difficult than the given proof (the geometric intersection number of
$C_{1}$ and $C_{2}$ is $72$).

\section{The proof of Theorem \ref{noextend}}

In this section we sketch the proof of Theorem \ref{noextend}
following \cite{JJ}. We include this sketch since \cite{JJ} did not
appear, and although examples as in Theorem \ref{noextend} appear to
be well-known, no explicit example appears to be recorded in the
literature.

We begin by considering a Heegaard embedding $h
: \Sigma_{g} \rightarrow S^{3}$ (ie, an embedding such that
$h(\Sigma_{g})$ bounds handlebodies on both sides).
We perturb this Heegaard splitting using an automorphism $f \in {\cal
I}(\Sigma_{g})$ to give a Heegaard splitting of a new manifold $M_{h,f}$
as follows.

We let $H_{+}$ and $H_{-}$ be the handlebodies on the positive and
negative sides of $h(\Sigma_{g})$ respectively (so the positive unit
normal
to $h(\Sigma_{g})$ points into $H_{+}$).
The manifold $M_{h,f}$ is constructed by gluing $\partial H_{-}$ to
$\partial H_{+}$ by
$$h \circ f \circ h^{-1} : \partial H_{-} \rightarrow  \partial H_{+}. $$
If $h \circ f \circ h^{-1}$ extended over the handlebody $H_{+}$, then one
can check that $M_{h,f} \cong S^{3}$.
In particular, suppose $f$ extends over {\em some} handlebody
$\widetilde{f} : H \rightarrow H$.
Then if $H_{+}$ and $H_{-}$ are the positive and negative handlebodies in
a genus $g$ Heegaard splitting of $S^{3}$, any diffeomorphism
$\widetilde{h} : H \rightarrow H_{+}$ restricts to a Heegaard embedding
$h: \Sigma_{g} \rightarrow S^{3}$, and $h \circ f \circ h^{-1}$ extends
over $H_{+}$ by $\widetilde{h} \circ \widetilde{f} \circ
\widetilde{h}^{-1}$.
So to prove that an automorphism $f$ does not extend over any handlebody,
it suffices to show that for any Heegaard embedding $h:\Sigma_{g}
\rightarrow S^{3}$, the manifold $M_{h,f}$ is not diffeomorphic to
$S^{3}$.

Next we note that since $f \in {\cal I}(\Sigma_{g})$, $M_{h,f}$ will
be an integral homology $3$-sphere.
We let $\mu(M_{h,f}) \in {\mathbb Z}/2{\mathbb Z}$ denote the Rohlin
invariant of $M_{h,f}$ (see eg \cite{GS}).
Since $\mu(S^{3}) = 0$, we wish to find $f \in {\cal I}(\Sigma_{g})$
such that for every Heegaard embedding $h$, $\mu(M_{h,f}) = 1$.  In
\cite{Jo}, Johnson studies the {\em Birman-Craggs homomorphisms} (see
\cite{BC}) and gives a very effective way of finding such $f$.

In what follows, the main ideas and relevant theorems of \cite{Jo}
necessary for our purpose are stated without proofs (see \cite{Jo} for the
proofs and complete references).

A {\em symplectic quadratic form} ($Sp$-form for short) on
$H_{1}(\Sigma_{g};{\mathbb Z}/2{\mathbb Z})$ is a function
$$\omega : H_{1}(\Sigma_{g};{\mathbb Z}/2{\mathbb Z}) \rightarrow {\mathbb
Z}/2{\mathbb Z}$$
such that
$$\omega(a+b) = \omega(a) + \omega(b) + a \cdot b$$
where $a \cdot b$ is the symplectic bilinear form given by the mod 2
intersection number of $a$ and $b$.
We denote the set of $Sp$-forms on $H_{1}(\Sigma_{g};{\mathbb Z}/2{\mathbb
Z})$ by $\Omega = \Omega(\Sigma_{g})$.

Given an element $a \in H_{1}(\Sigma_{g};{\mathbb Z}/2{\mathbb Z})$, we
obtain a function $\overline{a} : \Omega \rightarrow {\mathbb Z}/2{\mathbb
Z}$ defined by $\overline{a}(\omega)\! = \!\omega(a)$.
If $\{ a_{1},b_{1},...,a_{g},b_{g} \}$ is a symplectic basis for
$H_{1}(\Sigma_{g};{\mathbb Z}/2{\mathbb Z})$ (with respect to $\cdot$),
then the {\em Arf-invariant} of a form $\omega \in \Omega$ is defined to
be $Arf(\omega) = \left( \sum_{i=1}^{g}\overline{a_{i}}\overline{b_{i}}
\right) (\omega)$.
We denote the set of $Sp$-forms with zero Arf-invariant by $\Psi =
\Psi(\Sigma_{g}) \subset \Omega(\Sigma_{g})$.

Now let $h: \Sigma_{g} \rightarrow S^{3}$ be a Heegaard embedding.
Seifert's linking form defines an $Sp$-form by setting
$$\omega_{h}(a) = \lambda(h(\gamma_{a}),h(\gamma_{a})^{+})$$
where $\gamma_{a}$ is a simple closed curve representing the class $a \in
H_{1}(\Sigma_{g};{\mathbb Z}/2{\mathbb Z})$ and
$\lambda(h(\gamma_{a}),h(\gamma_{a})^{+})$ is the mod 2 linking number of
$h(\gamma_{a})$ and its push off in the positive normal direction,
$h(\gamma_{a})^{+}$.
The form $\omega_{h}$ lies in $\Psi$, and furthermore, any $\omega \in
\Psi$ can be realized by some Heegaard embedding of $\Sigma_{g}$.

The following facts from \cite{Jo} will essentially complete the proof
(labellings below are those of \cite{Jo}).
\begin{itemize}
\item {\bf Corollary 1 to Theorem 1}\qua$\mu(M_{h,f})$ depends only on
$\omega_{h}$ and $f$.
\item {\bf Lemma 11}\qua If we denote the Abelian group of functions from
$\Psi$ into ${\mathbb Z}/2{\mathbb Z}$ by ${\mathbb Z}/2{\mathbb
Z}^{\Psi}$, then there is a homomorphism
$$\sigma : {\cal I}(\Sigma_{g}) \rightarrow {\mathbb Z}/2{\mathbb
Z}^{\Psi}$$
such that if we denote the image of $f$ under $\sigma$ by $\sigma_{f}$,
then $\sigma_{f}(\omega_{h}) = \mu(M_{h,f})$.
\item {\bf Consequence of Theorem 4}\qua The constant function $1$ is in the
image of $\sigma$ for every $g \geq 2$.
\end{itemize}

It now follows that if we let $f \in {\cal I}(\Sigma_{g})$ be such that
$\sigma_{f} = 1$, then for any Heegaard embedding $h$, $\mu(M_{h,k}) = 1$,
and hence $f$ cannot extend over any handlebody.
This completes the proof of the first statement of Theorem \ref{noextend}

For the second, note the following consequence of the proof of the
first statement in Theorem \ref{noextend}, which gives the second
statement
of Theorem \ref{noextend}.

\begin{scholium} \label{its_a_scholium}
If $f \in {\cal I}(\Sigma_{g})$ satisfies $\sigma_{f} = 1$, then for any
odd integer $n$, $f^{n}$ does not extend over any handlebody.
\end{scholium}

\noindent
\proof Note that for any odd $n$,
$$\sigma_{f^{n}} = n \sigma_{f} = \sigma_{f} = 1$$
so $f^{n}$ cannot extend over any handlebody.\endproof

\subsection{} Here we extend Theorem \ref{noextend} to obtain
pseudo-Anosov maps for all genera $\geq 2$. In the
notation developed above, Theorem \ref{pA} follows from Scholium
\ref{its_a_scholium} and

\begin{theorem}
\label{pAnoextend}
For any $g \geq 2$, there exists pseudo-Anosov mapping classes $f \in
{\cal I}(\Sigma_{g})$ for which $\sigma_{f}=1$.\end{theorem}

\proof We use the notation of the proof of Theorem \ref{noextend}.
Let $h \in {\cal I}(\Sigma_{g})$ be such that $\sigma_{h} = 1$, and let
$\phi \in ker(\sigma)$ be a pseudo-Anosov mapping class with stable and
unstable laminations $[\lambda_{s}], [\lambda_{u}] \in {\cal
PML}(\Sigma_{g})$.
By choosing a different $h$ if necessary, we may assume that
$h([\lambda_{s}]) = [\mu] \neq [\lambda_{u}]$.

Let $n$ be a positive integer and $V_{s}, V_{u},$ and $V_{\mu}$ be
neighborhoods in ${\cal PML}(\Sigma_{g})$ of
$[\lambda_{s}],[\lambda_{u}],$ and $[\mu]$, respectively, so
that

\begin{itemize}
\item For all $[\nu] \in \overline{V}_{\mu}$ and $[\eta] \in
\overline{V}_{u}$, $\nu \pitchfork \eta$ and $\nu \cup \eta$ fills
$\Sigma_{g}$ (so, $\overline{V}_{\mu} \cap \overline{V}_{u} = \emptyset$).
\item $\phi^{n}({\cal PML}(\Sigma_{g}) \setminus V_{u}) \subset
\overline{V}_{s}$, $\phi^{-n}({\cal PML}(\Sigma_{g}) \setminus
V_{s}) \subset \overline{V}_{u}$, and $\overline{V}_{s} \cap
\overline{V}_{u} = \emptyset$
\item $h(\overline{V}_{s}) = \overline{V}_{\mu}$
\item $\overline{V}_{s} \cong \overline{V}_{u} \cong \overline{V}_{\mu}
\cong B^{n}$
\end{itemize}

The first property is possible to arrange since $\mu \pitchfork
\lambda_{u}$ and $\mu \cup \lambda_{u}$ fills $\Sigma_{g}$ implies the
existence of disjoint neighborhoods with the same property.
The second follows from standard properties of
the dynamical behavior of the action of $\phi$ on
${\cal PML}(\Sigma_{g})$ (see for example \cite{I}, Chapter 8).
The third is possible because $h([\lambda_{s}]) = [\mu]$ and $h$ acts by a
homeomorphism on ${\cal PML}(\Sigma_{g})$.
The last property is possible because ${\cal PML}(\Sigma_{g})$ is a
manifold.

We let $f = h\phi^{n} \in {\cal I}(\Sigma_{g})$ and note that
$$f({\cal PML}(\Sigma_{g}) \setminus V_{u}) = h(\phi^{n}({\cal
PML}(\Sigma_{g}) \setminus V_{u})) \subset h(\overline{V}_{s}) =
\overline{V}_{\mu}$$
and
$$f^{-1}({\cal PML}(\Sigma_{g}) \setminus V_{\mu}) =
\phi^{-n}(h^{-1}({\cal PML}(\Sigma_{g}) \setminus V_{\mu})) =
\phi^{-n}({\cal PML}(\Sigma_{g}) \setminus V_{s}) \subset
\overline{V}_{u}$$
In particular, $f(\overline{V}_{\mu}) \subset \overline{V}_{\mu}$ and
$f^{-1}(\overline{V}_{u}) \subset \overline{V}_{u}$.
By the Brouwer Fixed Point Theorem, $f$ has fixed points $[\nu_{s}] \in
\overline{V}_{\mu}$ and $[\nu_{u}] \in \overline{V}_{u}$.
Since $f$ is clearly not periodic, and because $\nu_{s} \pitchfork
\nu_{u}$ and $\nu_{s} \cup \nu_{u}$ fill $\Sigma_{g}$, it follows that $f$
must be pseudo-Anosov with stable and unstable laminations $[\nu_{s}]$ and
$[\nu_{u}]$ respectively.
Since $\sigma$ is a homomorphism, $\sigma_{f} = 1$.\endproof

\subsection{} Following $\cite{JJ}$, we construct some explicit 
examples of $f \in {\cal I}(\Sigma_{2})$ which do not extend over any
handlebody.

\begin{example} \label{bunchofem}
If $a,b \in H_{1}(\Sigma_{2};{\mathbb Z}/2{\mathbb Z})$ satisfy $a \cdot b
=1$, we can find a pair of transversely intersecting simple closed curves
$\alpha$ and $\beta$ representing $a$ and $b$ respectively such that
$\alpha \cap \beta$ is exactly $1$ point.
The regular neighborhood, $N(\alpha \cup \beta)$ is homeomorphic to a
torus-minus-disk embedded in $\Sigma_{2}$ and $\gamma = \partial N(\alpha
\cup \beta)$ is a separating essential simple closed curve.
If we denote a Dehn twist about $\gamma$ by $T_{\gamma}$ (note that
$T_{\gamma} \in {\cal I}(\Sigma_{2})$) then according to \cite{Jo} (Lemma
12a)
$$\sigma_{T_{\gamma}} = \overline{a} \overline{b}.$$
Now fix a symplectic basis $a_{1},b_{1},a_{2},b_{2}$ for
$H_{1}(\Sigma_{2};{\mathbb Z}/2{\mathbb Z})$, and let
$\gamma_{1},...,\gamma_{10}$ be separating essential simple closed curves
associated, as above, to the following $10$ pairs of elements $a,b \in
H_{1}(\Sigma_{2};{\mathbb Z}/2{\mathbb Z})$ (with $a \cdot b = 1$).
$$\begin{array}{lll}
a_{1},b_{1} \leadsto \gamma_{1} \, , & 
a_{1},b_{1}+a_{2} \leadsto \gamma_{2} \, , &
a_{1},b_{1}+b_{2} \leadsto \gamma_{3} \, , \\
b_{1},a_{1}+a_{2} \leadsto \gamma_{4} \, , &
a_{1},b_{1}+a_{2}+b_{2} \leadsto \gamma_{5} \, , &
b_{1},a_{1}+b_{2} \leadsto \gamma_{6} \, , \\
b_{1},a_{1}+a_{2}+b_{2} \leadsto \gamma_{7} \, , &
a_{1}+b_{1},a_{1}+a_{2} \leadsto \gamma_{8} \, , &
a_{1}+b_{1},a_{1}+b_{2} \leadsto \gamma_{9} \, , \\
a_{1}+b_{1},a_{1}+a_{2}+b_{2} \leadsto \gamma_{10}. & & \\ \end{array}$$
To compute $\sigma_{T_{\gamma_{j}}}$ for each $j=1,...,10$, we note first
that the defining characteristic of $Sp$-forms implies $\overline{a + b} =
\overline{a} + \overline{b} + a \cdot b$.
We also note that if $\phi: \Psi \rightarrow {\mathbb Z}/2{\mathbb Z}$,
then $\phi^{2} = \phi$.
It then follows that
$$\begin{array}{ll}
\sigma_{T_{\gamma_{1}}} = \overline{a}_{1} \overline{b}_{1} \, , &
\sigma_{T_{\gamma_{2}}} = \overline{a}_{1}\overline{b}_{1} + \overline{a}_{1}\overline{a}_{2} \, , \\
\sigma_{T_{\gamma_{3}}} = \overline{a}_{1}\overline{b}_{1} + \overline{a}_{1}\overline{b}_{2} \, , &
\sigma_{T_{\gamma_{4}}} = \overline{b}_{1}\overline{a}_{1} + \overline{b}_{1}\overline{a}_{2} \, , \\
\sigma_{T_{\gamma_{5}}} = \overline{a}_{1}\overline{b}_{1} + \overline{a}_{1}\overline{a}_{2} + \overline{a}_{1}\overline{b}_{2} + \overline{a}_{1} \, , &
\sigma_{T_{\gamma_{6}}} = \overline{b}_{1}\overline{a}_{1} + \overline{b}_{1}\overline{b}_{2} \, , \\
\sigma_{T_{\gamma_{7}}} = \overline{b}_{1}\overline{a}_{1} + \overline{b}_{1}\overline{a}_{2} + \overline{b}_{1}\overline{b}_{2} + \overline{b}_{1} \, , &
\sigma_{T_{\gamma_{8}}} = \overline{a}_{1}\overline{b}_{1} + \overline{a}_{1}\overline{a}_{2} + \overline{b}_{1}\overline{a}_{2} + \overline{a}_{2} \, , \\
\sigma_{T_{\gamma_{9}}} = \overline{a}_{1}\overline{b_{1}} + \overline{a}_{1}\overline{b}_{2} + \overline{b}_{1}\overline{b}_{2} + \overline{b}_{2} \, , &
\sigma_{T_{\gamma_{10}}} = \overline{a}_{1}\overline{b}_{1} + \overline{a}_{1}\overline{a}_{2} + \overline{a}_{1}\overline{b}_{2} + \overline{b}_{1}\overline{a}_{2} + \overline{b}_{1}\overline{b}_{2} \\
 & + \overline{a}_{1} + \overline{b}_{1} + \overline{a}_{2} + \overline{b}_{2} + 1.\ \end{array}$$
Any word in $T_{\gamma_{1}},...,T_{\gamma_{10}}$ such that the total
exponent of each $T_{\gamma_{j}}$ is odd provides an automorphism in
${\cal I}(\Sigma_{2})$ which does not extend over any handlebody.
\end{example}

\subsection{Proof of Corollary \ref{newblob}}

Let $\theta$ be any word in $T_{\gamma_{1}},...,T_{\gamma_{10}}$ such
that the total exponent of each $T_{\gamma_{j}}$ is odd. Let 
$n\geq 1 $ be an integer, define $\theta_n =  T_{\gamma_1}^{2n}\theta$, and
let $M_n$ denote the mapping torus of $\theta_n$.
By the remarks in \S 6.2, $\theta_n$ does not extend over any handlebody.
Note that by the Remark at the end of \S 3, $\theta_n$ is pseudo-Anosov,
however we require the following description to gain extra control
of commensurability. By Lemma 1.1 of \cite{LM} the manifolds $M_n$
can be described as surgeries on a 1 cusped hyperbolic 3-manifold. Since
the degree of the invariant trace-field gets arbitrarily large on such
a sequence of surgeries (see \cite{LR}) and the invariant trace-field
is an invariant of the commensurability class, by subsequencing if necessary
we obtain the set of non-commensurable manifolds.\endproof

\Addresses
\end{document}